\documentclass[a4paper,leqno]{article}
\usepackage{latexsym,amssymb,amsmath}
\usepackage{color,amsfonts}
\usepackage[english]{babel}

\newtheorem{theorem}{Theorem}[section]
\newtheorem{lemma}[theorem]{Lemma}
\newtheorem{proposition}[theorem]{Proposition}
\newtheorem{corollary}[theorem]{Corollary}

\newtheorem{remark}[theorem]{Remark}
\newenvironment{proof}[1][Proof]{\textbf{#1.} }{\ \mathbb Rule{0.5em}{0.5em}}

\begin{document}

\title{Inviscid limit of stochastic damped 2D Navier-Stokes equations}
\author{Hakima Bessaih\footnote{University of Wyoming, Department of Mathematics, Dept. 3036, 1000
East University Avenue, Laramie WY 82071, United States, bessaih@uwyo.edu} \ 
\& Benedetta Ferrario\footnote{Universit\`a di Pavia, Dipartimento 
di Matematica, via  Ferrata 1, 27100 Pavia, Italy, benedetta.ferrario@unipv.it}}

\date{}
\maketitle

\begin{abstract}
We consider the inviscid limit of the stochastic damped 
2D Navier-Stokes equations. 
We prove that, when the viscosity vanishes, 
the stationary solution of the stochastic 
damped Navier-Stokes equations converges to a stationary solution
of the stochastic damped Euler equation 
and that  the rate of dissipation 
of enstrophy converges to zero.
In particular, this limit obeys an enstrophy balance.
The rates are computed with respect to a limit measure of the unique invariant measure of
the stochastic damped Navier-Stokes equations.
\end{abstract}
\noindent
{\bf MSC2010}: 60G10, 60H30, 35Q35.\\
{\bf Keywords}: Inviscid limits, enstrophy balance, 
stationary processes, invariant measures.

\section{Introduction} 
In this paper, we are interested in the equations of motion of 
incompressible fluids in a bounded domain of $\mathbb R^2$. 
In particular, we consider the Euler or Navier-Stokes equations
damped by a term proportional to the velocity. 
Damping terms in two dimensional turbulence studies have been
considered to model pumping due to friction with boundaries. 
Numerical studies of two dimensional turbulence employ devices to remove
the energy that piles up at the large scales, and damping is the 
most common such device.
We refer to  \cite{Gal, boffetta} for a 
physical motivation of the model and to 
\cite{BCT88, ilyin1, ilyin2} for a mathematical analysis of the deterministic 
damped Navier-Stokes equations and to \cite{Bes2000, Bes2008} 
for the stochastic  damped Euler equations. 
 
These stochastic damped equations are given by
\begin{equation}\label{stoc-damp}
\begin{cases}
du+[-\nu \Delta u+(u\cdot \nabla)u +\gamma u +\nabla p]dt=dw&\\
\nabla \cdot u=0&
\end{cases}
\end{equation}
The non negative
coefficients $\nu$ and $ \gamma$ are called kinematic viscosity and sticky
viscosity, respectively.  The unknowns are the velocity $u$ and the
pressure $p$.
Suitable boundary conditions have to be considered; 
in this paper the spatial domain is a box and 
periodic boundary conditions are assumed.

For a  fixed  $\gamma>0$, if $\nu>0$ these are called the stochastic 
damped Navier-Stokes equations, whereas if $\nu=0$
they are the stochastic damped Euler equations. 
If $\gamma=0$ and $\nu=0$, we refer to
\cite{Bes,Pez,cutland, CFM, GHV, kim, Mikulevic} 
for an analysis of the existence and/or uniqueness of solutions
and 
to \cite{eyink} where some dissipation of enstrophy 
arguments are discussed in Besov spaces.

Turbulence theory investigates the  behavior of certain quantities
as the viscosity $\nu$ vanishes. 
In particular, in the two dimensional setting one is interested 
in understanding what happens to 
the balance equation of energy and enstrophy (in the stationary
regime) as the viscosity vanishes.
D. Bernard \cite{bernard} suggested that 
there is no anomalous dissipation of enstrophy
in damped and driven Navier-Stokes equations;
Constantin and Ramos \cite{CR} proved that there 
is no anomalous dissipation neither of energy nor of 
enstrophy as $\nu \to 0$
for the deterministic damped Navier-Stokes equations in the whole 
plane. Some similar
questions were suggested by Kupiainen  \cite{K} for the stochastic case.
Therefore we address the same problem when the forcing  term is of  
white noise type. Tools from stochastic analysis 
are very useful to investigate the 
same problem studied in \cite{CR}, giving a rigorous 
meaning to the averages of velocity and vorticity.
Indeed, using stochastic PDE's allows 
to express the stationary regime by means of an invariant measure,
whereas in the deterministic setting the stationary regime is 
described by taking time averages on the infinite time interval.

In this paper we shall prove that in the stationary regime
system \eqref{stoc-damp} has 
no anomalous dissipation neither of energy nor of enstrophy 
as $\nu \to 0$. 
However, we shall be working in a finite two dimensional spatial 
domain and not in the whole plane; this answers 
one of the questions posed by Kupiainen
in \cite{K} about the behaviour of the stochastic damped Navier-Stokes
equations on a torus for vanishing viscosity.

As far as the content of the paper is concerned, in Section 2 we introduce 
some functional spaces,
the equations in their vorticity formulation and the assumptions 
on the noise term.  We also introduce 
the classical properties of the nonlinear term associated to these 
equations.  Section 3 is devoted to the well posedness of the 
stochastic 2D damped Navier-Stokes equations, where some uniform estimates 
are computed.  Starting from a known result of  existence and uniqueness 
of the invariant measure, we provide a balance law for the enstrophy.
The vanishing viscosity limit is studied in Section 4 and stationary 
solutions are constructed by means of a tightness argument providing 
a balance relation for these stationary solutions. Using these results, 
we provide a proof of no anomalous of enstrophy and energy for the 
stochastic damped 2D Navier-Stokes equations.

\section{Notations and  hypothesis}
Let the spatial domain $D$ be the square $[-\pi,\pi]^2$; periodic boundary
conditions are assumed. 
A basis of the space $L^2(D)$ with periodic boundary
conditions is $\{e_k\}_{k \in \mathbb Z^2}$, 
$e_k(x)=\frac 1{2\pi}e^{i k\cdot x}$, whereas a basis for the space of
periodic vector fields  which are square integrable and divergence free is 
$\{\frac{k^\perp}{|k|} e_k\}_{k \in \mathbb Z^2}$, being $k^\perp=(-k_2,k_1)$.
Actually we consider $k \neq (0,0)$, since if $u$ is a solution of
system \eqref{stoc-damp} then also $u+c$ is a solution for any $c\in
\mathbb R$. Therefore we consider velocity fields with vanishing  mean value.

Let $\mathbb Z^2_0=\mathbb Z^2\setminus\{(0,0)\}$, and
$\mathbb Z^2_+=\{k=(k_1,k_2)\in \mathbb Z: k_1>0\}
\cup \{k=(0,k_2)\in \mathbb Z^2: k_2>0\}$.
Given $x=(x_1,x_2) \in \mathbb R^2$ we denote by $|x|$ its  norm: 
$|x|=\sqrt{(x_1)^2+(x_2)^2}$. Given $y=\Re y +i \Im y \in \mathbb C$ we denote  
by $|y|$ its absolute value and by $\overline y$ its complex
conjugate: $|y|=\sqrt{(\Re y)^2+(\Im y)^2}$, $\overline y=\Re y-i \Im y$.

For any $a \in \mathbb R$ we define the Hilbert space
\[
 H^a=\{f=\sum_{k \in \mathbb Z^2_0} f_k e_k(x): 
 \sum_{k \in \mathbb Z^2_0} |f_k|^2|k|^{2a}<\infty\}
\]
with scalar product 
$$
\langle f,g\rangle_{H^a}=\sum_{k \in \mathbb Z^2_0} |k|^{2a}f_k \overline {g_k};
$$
we set
$$
 \|f\|^2_{H^a}=\sum_{k \in \mathbb Z^2_0} |k|^{2a}|f_k|^2.
$$
For a vector $f=(f_1,f_2)$ we set
$$
 \|f\|^2_{H^a}=\|f_1\|^2_{H^a}+\|f_2\|^2_{H^a}.
$$
In particular, for scalar functions we have 
$\|f\|^2_{H^0}=\|f\|^2_{L^2(D)}$ and $\|f\|^2_{H^1}=\|\nabla f\|^2_{H^0}$.

The space $H^a$ is compactly embedded in the space $H^b$ if $a>b$.

Moreover, we consider the Banach spaces $W^{1,q}(D)$ ($1 \le q \le\infty$) endowed with the norm
$$
 \|f \|^q_{W^{1,q}(D)}= \|f \|^q_{L^q}+ \|\nabla f \|^q_{L^q}
$$
where $\|\cdot \|^q_{L^q}$ is the $L^q(D)$-norm.

Given a separable Hilbert space $X$, for $\alpha>0$ and $p\ge 1$
 we define the Banach space
$$
 W^{\alpha,p}(0,T;X)= 
 \left\{f \in L^p(0,T; X):\int_0^T \int_0^T 
   \frac{\|f(t)-f(s)\|_X^p}{|t-s|^{1+p\alpha}} dt\ ds <\infty \right\}
$$
and we set 
$$
 \|f\|^p_{W^{\alpha,p}(0,T;X)}= \int_0^T \|f(t)\|^p_X dt+
  \int_0^T \int_0^T    \frac{\|f(t)-f(s)\|_X^p}{|t-s|^{1+p\alpha}} dt\ ds .
$$

Let $(\Omega,F,P)$  be a complete probability space, 
with expectation denoted by $\mathbb E$.
We assume that the stochastic forcing term in \eqref{stoc-damp} is of the form
\[
 w=w(t,x)
 =
 \sum_{k \in \mathbb Z^2_0}\sqrt {q_k}\beta_k(t) \frac {k^{\perp}}{|k|}e_k(x).
\]
Here $\left\{\beta_{k}\right\}_{k\in \mathbb Z^2_+}$ is a sequence 
of independent complex-valued  standard 
Brownian motions on  $(\Omega,F,P)$, i.e. 
$\beta_{k}(t)=\Re\beta_k(t)+i\Im\beta_k(t)$ with $\{\Re\beta_k\}\cup
\{\Im\beta_k\}_{k \in \mathbb Z^2_+}$ a sequence of 
independent standard real Brownian motions; moreover we set
$\beta_{-k}=-\overline{ \beta_k}$ and $q_k=q_{-k}$ for any 
$k \in \mathbb Z^2_+$. Therefore  
$$w(t,x)=2\sum_{k \in \mathbb Z^2_+}
\sqrt {q_k} \frac {k^{\perp}}{|k|} \left[\Re\beta_k(t)\cos(k \cdot x) 
- \Im\beta_k(t) \sin(k \cdot x)\right].
$$

In the 2D setting it is convenient to introduce the (scalar) vorticity
\[
 \xi=\nabla^\perp \cdot  u\equiv \frac{\partial u_2}{\partial x_1} -
\frac{\partial u_1}{\partial x_2}.
\]
System \eqref{stoc-damp} corresponds to 
\begin{equation}\label{eq-vort}
\begin{cases}
d\xi+[-\nu \Delta \xi+\gamma \xi +u\cdot \nabla\xi ] dt=dw^{curl}&\\
\xi=\nabla^\perp \cdot u&
\end{cases}
\end{equation}
obtained by taking the curl of both sides of the first equation of \eqref{stoc-damp}.
Periodic boundary conditions have to be added to this system.
The noise is $w^{curl}(t,x)=-2 \sum_{k \in \mathbb Z^2_+} \sqrt {q_k} |k| 
\left[\Im\beta_k(t)\cos(k \cdot x)+\Re\beta_k(t)\sin(k\cdot x)\right]$.
Let us define 
\begin{equation}\label{noise2}
 Q:=\sum_{k\in \mathbb Z^2_0} |k|^2 q_k .
\end{equation}  

Classical results are
\begin{equation}
 \mathbb E\|w^{curl}(t)\|^2_{H^0}=2t Q \qquad \forall t\ge 0
\end{equation}
\begin{equation}\label{wcurl}
 \mathbb E\|w^{curl}\|^p_{W^{\alpha,p}(0,T;H^0)}\le C(\alpha,p) (T^{1+p/2}+1) (Q)^{p/2}
\end{equation}
for any $\alpha\in (0,\frac 12)$, $p\ge 2$, and
the Burkh\"older-Davies-Gundy inequality
\begin{equation}\label{BDG}
 \mathbb E \left(\sup_{0\le t\le T}\int_0^t
 \langle |\xi(s)|^{p-2}\xi(s),dw^{curl}(s)\rangle_{L^{2}}\right)
 \le C(p) \sqrt Q \ \mathbb E \sqrt{ \int_0^T \|\xi(s)\|^{2(p-1)}_{L^p}ds } 
\end{equation}
For this latter inequality we have used that $\sup_{x \in
  D}|e_k(x)|=1$ for all $k$.
\\
Here and henceforth,
$C(\cdot)$ denotes a positive constant depending on the specified 
parameters; it may change from line to line. 

Knowing the vorticity $\xi$, we recover the velocity $u$ by solving the
elliptic equation 
\begin{equation}\label{ell}
 -\Delta u = \nabla^{\bot}\xi.
\end{equation}
This means that if $\xi(x)=\sum_k \xi_k e_k(x)$, then 
$u(x)=-i\sum_k\frac{k^\perp}{|k|^2} \xi_k e_k(x) $.

We present basic properties of the bilinear term $u \cdot \xi$ in the 2D
setting. These are classical results in the analysis of incompressible
fluids (see e.g. \cite{Tem}).
\begin{lemma}
There exists a positive constant $C$ such that
\begin{equation}\label{stima-tril-u}
\Big|\int_D (u \cdot \nabla) v \cdot \psi\ dx \Big|
\le 
C \|u\|_{L^4} \|v\|_{L^4} \|\psi\|_{H^1}
\end{equation}
for all  divergence free vectors with the regularity specified in the r.h.s.,
and for any $a>1$ 
\begin{equation}\label{stima-tril-xi}
 \Big|\int_D u \cdot \nabla \xi \ \phi \ dx\Big| \le C 
 \|u\|_{H^0} \|\xi\|_{H^1}\|\phi\|_{H^a},
\end{equation}
\begin{equation}\label{stima-tril-xi0}
\Big|\int_D u \cdot \nabla \xi \ \phi \ dx\Big| \le C 
 \|u\|_{H^0} \|\xi\|_{H^{1+a}}\|\phi\|_{H^0}   
\end{equation}
for all  functions with the regularity specified in the r.h.s..
\end{lemma}
\proof The key relationship for \eqref{stima-tril-u} is
$$
 \int_D [u \cdot \nabla]v \cdot \psi\ dx =-
  \int_D [u \cdot \nabla]\psi \cdot v\ dx 
$$
assuming sufficient regularity for $u,v,\psi$; 
this is obtained by integrating by parts. 
Then, we get the estimate by H\"older inequality and this
is extended by density to vectors with the specified regularity.
For \eqref{stima-tril-xi} we use  H\"older inequality
and the continuous embedding  $H^{a}\subset L^\infty(D)$ for $a>1$.
Similarly, we obtain the latter estimate.
\hfill $\Box$

\begin{lemma}\label{enstr}
Let $\xi=\nabla^\perp \cdot  u$. We have
\begin{equation}\label{cambiasegno}
 \int_D [u\cdot \nabla \xi] \phi\ dx=-
 \int_D [u\cdot \nabla \phi] \xi\ dx  
\qquad \forall \xi,\phi \in H^1
\end{equation}
and for any $p \ge 2$
\begin{equation}\label{valori-distribuz-q}
 p \int_D [u\cdot \nabla \xi]\ \xi |\xi|^{p-2} \psi\ dx=-
\int_D [u \cdot\nabla \psi]\ |\xi|^p \ dx
\qquad \forall \xi \in L^{2p}, \psi \in H^1 .
\end{equation}
Moreover,
\begin{equation}\label{cons}
\int_D [u\cdot \nabla \xi] \xi\ |\xi|^{p-2}  dx=0
\qquad \forall \xi \in L^{2p} .
\end{equation}

\end{lemma}
\proof The two first relationships  
\eqref{cambiasegno}-\eqref{valori-distribuz-q} are easily 
obtained by integrating by parts, where in \eqref{valori-distribuz-q} the proof
is done first with smooth functions and then by density it is extended
on the spaces specified; notice that for $p>1$, if $\xi \in L^{2p}$ 
then $u \in W^{1,2p}\subset L^\infty$ 
and the r.h.s. is meaningful (see \cite{LML2006} ). 
Eventually, \eqref{cons} is the
particular case of \eqref{valori-distribuz-q} for 
$\psi=1$. \hfill $\Box$

\section{The stochastic damped Navier-Stokes equations}
The well posedness of the stochastic damped 2D Navier-Stokes equations 
\begin{equation}\label{DNSE}
\begin{cases}
 d\xi^\nu+[-\nu \Delta \xi^\nu+u^\nu\cdot \nabla\xi^\nu +\gamma\xi^\nu] dt
    =dw^{curl}&\\
 \xi^\nu=\nabla^\perp \cdot u^\nu&
\end{cases}
\end{equation}
is very similar to the case when $\gamma=0$. Here, we assume periodic
boundary conditions with period box $[-\pi,\pi]^2$.

The proof of  existence of a unique solution 
for square summable initial vorticity is the same as
 the proof for square summable initial velocity that
can be found in \cite{Fladiss}, where the proof is performed  for $\gamma=0$. 
Similar proofs can also be found in \cite{Bes, Pez} 
with some uniform estimates with respect to the viscosity $\nu$.
Here, we point out the peculiar estimate \eqref{stima22} for $\gamma>0$,  
useful in the analysis of the limit as $\nu \to 0$.

\begin{theorem}\label{theorem1}
Let $\gamma, \nu>0$, $p\ge 2$. 
Assume 
$$
\mathbb E \|\xi^\nu(0)\|_{L^p}^p<\infty,
\qquad
Q<\infty. 
$$ 
Then, there exists  
a process $\xi^\nu$ with paths in 
$C([0,\infty), L^p)\cap L^{2}_{loc}(0,\infty; H^1)$ $P$-a.s., 
which is a Feller Markov process in $L^p$ 
and is the unique solution for \eqref{DNSE} with initial data $\xi^\nu(0)$. 
Moreover,  there exist two positive constants $C(p,T)$ and $C(p)$,
independent of $\nu$, such that
\begin{equation}\label{stima-nu}
\mathbb E \sup_{0\le t\le T}\|\xi^\nu(t)\|^p_{L^p}\le C(p,T)
\end{equation}
for any finite $T$, and
\begin{equation}\label{stima22}
 \sup_{0\le t<\infty}\mathbb E \|\xi^\nu(t)\|^p_{L^p}\le C(p).
\end{equation}
In particular,  the constants  depend also on
 $\gamma, Q, \mathbb E\|\xi^\nu(0)\|_{L^p}^p$.
\end{theorem}
\proof 
The proof of the existence of solutions, which is quite classical
requires some Galerkin 
approximation of $\xi^{\nu}$, say  $\xi^{\nu,n}$, for which a priori  
estimates are proved uniformly in $n$. Using a subsequence of 
$\xi^{\nu,n}$ which converges in the weak  or
weak-star topologies of appropriate spaces, 
one can then prove that there exists a solution to 
\eqref{DNSE}. The proof of uniqueness and Feller property
is standard and hence  omitted.

Let $\nu>0$, $x \in D$ and $t\in [0,T]$; It\^o formula 
for $|\xi^\nu(t,x)|^p$ gives
\[
d|\xi^\nu(t,x)|^p = p |\xi^\nu(t,x)|^{p-2} \xi^\nu(t,x) d\xi^\nu(t,x) + \frac 12 p
(p-1)|\xi^\nu(t,x)|^{p-2}  2Q dt
\]
hence
\[
d|\xi^\nu(t,x)|^p + p |\xi^\nu(t,x)|^{p-2} \xi^\nu(t,x) 
[-\nu \Delta \xi^\nu(t,x)+u^\nu \cdot \nabla \xi^\nu(t,x)
+\gamma \xi^\nu(t,x)]\ dt 
\]
\[
-  p  (p-1)
|\xi^\nu(t,x)|^{p-2}  Q dt=p |\xi^\nu(t,x)|^{p-2} \xi^\nu(t,x) dw^{curl}(t,x)
\]
Integrating on the spatial domain $D$, by using \eqref{cons} and by integrating by parts we get
\begin{multline}\label{Lp-formula}
d\|\xi^\nu(t)\|_{L^p}^p+p\nu (p-1)\|\ |\xi^\nu(t)|^{\frac{p-2}2} \nabla \xi^\nu(t)\|^2_{H^0}dt
+p \gamma \|\xi^\nu(t)\|_{L^p}^p dt 
\\
-  Q p(p-1) \|\xi^\nu(t,x)\|_{L^{p-2}}^{p-2}   dt=
p \langle |\xi^\nu(t)|^{p-2} \xi^\nu(t),dw^{curl}(t)\rangle .
\end{multline}
Integrating over the finite time interval $(0,s)$
we get that
\begin{multline}\label{primastima}
 \|\xi^\nu(s)\|^p_{L^p}+\nu p(p-1)\int_{0}^{s}\| \ |\xi^\nu(r)|^{\frac{p-2}2}   \nabla\xi^\nu(r)\|^{2}_{H^{0}}dr
  +\gamma p\int_{0}^{s}\|\xi^\nu(r)\|^{p}_{L^p}dr\\
  = \|\xi^\nu(0)\|^p_{L^p}
 + p\int_{0}^{s} \langle  |\xi^\nu(r)|^{p-2} \xi^\nu(r), dw^{curl}(r)\rangle\\
  +Q p(p-1)\int_{0}^{s}\|\xi^\nu(r)\|^{p-2}_{L^{p-2}} dr.
\end{multline}
Therefore
\begin{multline}\label{nuova-unif}
 \sup_{0\le s \le T}\|\xi^\nu(s)\|^p_{L^p}
 \le \|\xi^\nu(0)\|^p_{L^p}
 + p\sup_{0\le s \le T}\int_0^s \langle  |\xi^\nu(r)|^{p-2} \xi^\nu(r), dw^{curl}(r)\rangle\\
 + Q p(p-1)\int_0^T \sup_{0\le r \le s}\|\xi^\nu(r)\|^{p-2}_{L^{p-2}} ds.
\end{multline}
On the other side, using first  Burkholder-Davis-Gundy inequality \eqref{BDG} and 
then  H\"older inequality, we have that

\begin{equation*}
\begin{split}
p\mathbb E\Big( \sup_{0\le s\le T}& \int_{0}^{s}\langle  |\xi^\nu(r)|^{p-2} \xi^\nu(r), dw^{curl}(r)\rangle\Big)\\
&\le 
 pC(p) \sqrt Q  \mathbb E\sqrt{ \int_0^T \|\xi^\nu(r)\|^{2p-2}_{L^p}  dr}\\
&\le p C(p) \sqrt Q  \mathbb E\left(\sup_{0\le s\le T}  \|\xi^\nu(s)\|^{p/2}_{L^p}
 \sqrt{  \int_{0}^{T} \|\xi^\nu(r)\|^{p-2}_{L^p}dr}\right)\\
&\le \frac{1}{2}\mathbb E\sup_{0\le s\le T} \|\xi^\nu(s)\|^{p}_{L^p}
 +\frac{Q}{2}C(p)^2 p^2\mathbb E
 \int_{0}^{T}\|\xi^\nu(r)\|^{p-2}_{L^p}dr\\
& \le \frac{1}{2}\mathbb E\sup_{0\le s\le T} \|\xi^\nu(s)\|^{p}_{L^p}
 +\frac{Q}{2} C(p)^2 p^2\mathbb E
 \int_{0}^{T}\sup_{0\le r \le s}\|\xi^\nu(r)\|^{p-2}_{L^p}ds.
\end{split}
\end{equation*}
Taking expectation in \eqref{nuova-unif} and collecting all the estimates
we get
\begin{equation}\label{E1}
\begin{split}
 \frac{1}{2}\mathbb E\sup_{0\le s\le T} &\|\xi^\nu(s)\|^p_{L^p}
 \le \mathbb E \|\xi^\nu(0)\|^p_{L^p}
  +Q C(p)  
  \int_{0}^{T}\mathbb E\sup_{0\le r\le s}\|\xi^\nu(r)\|^{p-2}_{L^p}ds\\
 & \le \mathbb E \|\xi^\nu(0)\|^p_{L^p}
  +\epsilon \int_0^T \mathbb E \sup_{0\le r\le s}\|\xi^\nu(r)\|^{p}_{L^p}ds
  +C(\epsilon,p, Q)T
\end{split}  
\end{equation}
for any $\epsilon>0$, by Young inequality. Using Gronwall lemma 
we obtain \eqref{stima-nu}. Taking expectation in \eqref{primastima}
and using \eqref{stima-nu}, we also get that
\begin{equation*}
 \nu (p-1)\mathbb E\int_{0}^{T}\| \ |\xi^\nu(s)|^{\frac{p-2}2} \nabla \xi^\nu(s)\|^{2}_{H^{0}}ds
 +\gamma \mathbb E\int_{0}^{T}\|\xi^\nu(s)\|^{p}_{L^p}ds
 \le C\left(p, T, Q, E \|\xi^\nu(0)\|^p_{L^p}\right).
\end{equation*}
For $p=2$ this gives in particular
\[
\mathbb E\int_{0}^{T} \|\nabla \xi^\nu(s)\|^{2}_{H^{0}}ds\le 
C\left(T, Q, E \|\xi^\nu(0)\|^2_{L^2}\right).
\]

Going back to estimate \eqref{primastima} and taking expectation, we
have
\begin{equation}
\begin{split}
\mathbb E \|\xi^\nu(s)\|^p_{L^p}&
 +\gamma p\int_{0}^{s}\mathbb E \|\xi^\nu(r)\|^{p}_{L^p}dr\\
& \le \mathbb E \|\xi^\nu(0)\|^p_{L^p}
 + Q p(p-1)\int_{0}^{s}\mathbb E \|\xi^\nu(r)\|^{p-2}_{L^{p-2}} dr\\
& \le  \mathbb E \|\xi^\nu(0)\|^p_{L^p} 
 + \frac {\gamma p }2 \int_{0}^{s}\mathbb E \|\xi^\nu(r)\|^{p}_{L^p} dr
 + C(\gamma, p , Q)s.
\end{split}  
\end{equation}
Hence
$$
\mathbb E \|\xi^\nu(s)\|^p_{L^p} 
 \le \mathbb E \|\xi^\nu(0)\|^p_{L^p} 
 -\frac {\gamma p }2 \int_{0}^{s}\mathbb E \|\xi^\nu(r)\|^{p}_{L^p} dr
 + C(\gamma, p , Q)s;
$$
Gronwall lemma gives 
$$
 \mathbb E \|\xi^\nu(s)\|^p_{L^p} \le \mathbb E
 \|\xi^\nu(0)\|^p_{L^p}e^{-\gamma ps/2 }+\frac{2C(\gamma,p,Q)}{\gamma p}
 \left(1-e^{-\gamma ps/2}\right)
$$
for any $s \in [0,\infty)$.
This implies \eqref{stima22}.
 \hfill $\Box$

\begin{remark}
The solution $\xi^\nu$ is a process whose paths are a.s. in 
$C([0,\infty), H^0)\cap L^{2}_{loc}(0,\infty; H^1)$ at least; 
therefore it solves 
system \eqref{DNSE} in the following
sense: for all $t \in [0,\infty)$ and $\phi \in H^a$ with $a>1$, we have
\begin{multline*}
 \int_D\xi^\nu(t,x)\phi(x)\  dx
 +\nu \int_0^t \int_D \nabla\xi^\nu(s,x)\cdot \nabla \phi(x)\ dx\  ds
 \\+\int_0^t \int_D u^\nu(s,x)\cdot \nabla\xi^\nu(s,x)\phi(x)\ dx \  ds
 +\gamma \int_0^t\int_D \xi^\nu(s,x) \phi(x)\ dx \ ds
 \\=\int_D \xi^\nu(0,x)\phi(x)\ dx+\int_D w^{curl}(t,x)\phi(x)\ dx
\qquad P-a.s.
\end{multline*}
The trilinear term 
is well defined thanks to \eqref{ell} and \eqref{stima-tril-xi}.

Moreover, let us denote by 
$\xi^\nu(\cdot;\eta)$ the solution with initial data $\eta$ 
and by $B_b(L^p)$, $C_b(L^p)$ the spaces of Borel
bounded functions, respectively  continuous and bounded functions, 
$\phi: L^p \to \mathbb R$.
To say that the  solution is a Feller process in $L^p$ (the $p$ depends
on the assumption on the initial vorticity) means that the Markov
semigroup
$P^\nu_t: B_b(L^p)\to B_b(L^p)$, defined as
\[
 \left(P^\nu_t\phi\right)(\eta)=\mathbb E \left[\phi(\xi^\nu(t;\eta))\right],
\]
actually maps $C_b(L^p)$ into itself.

We finally recall what is an invariant measure $\mu^\nu$: 
\[
 \int P_t^\nu \phi \ d\mu^\nu= \int  \phi \ d\mu^\nu \qquad \forall t \ge
 0, \phi \in L^p.
\]
The Feller property is important to prove the existence of invariant
measures by means of Krylov-Bogoliubov method (see, e.g., \cite{DPZ}).

\end{remark}

For any $\gamma> 0$ one can prove existence and uniqueness of the
invariant measure for system \eqref{DNSE}, 
following the lines of the proofs for the 2D
Navier-Stokes equation (the case $\gamma=0$). 
Indeed, Krylov-Bogoliubov method provides a way to prove the
existence of an invariant measure; this applies for a wide class of noises.
On the other side, uniqueness is a more delicate question. We just  recall
the best result of 
uniqueness of the invariant measure, proved by Hairer and
Mattingly \cite{HM06}. They assume
that the noise acts on first few modes, i.e.
\begin{equation}\label{Qf}
\begin{cases}
 \exists \mathcal Z \text{ finite } : 
 q_k \neq 0 \qquad\forall k \in \mathcal Z, \qquad\qquad q_k =0 \qquad\forall k \notin \mathcal Z&
\\
\text{ where $\mathcal Z$ has to be chosen in such a way that }&
\\
\quad \bullet
 \text{ it contains  at least two elements with different norms}&
\\
 \quad \bullet
 \text{ the integer linear combinations of elements of } \mathcal Z 
  \text{ generates }\mathbb Z^2&
\end{cases}
\end{equation}
Actually the kind and the number of forced modes, i.e. the elements of
$\mathcal Z$, is chosen independently of the viscosity.

We summarize the result.
\begin{theorem}\label{teo-mu-unica}
Let $\gamma>0$ and $2\le p <\infty$.
If \eqref{Qf} holds, 
then for any $\nu>0$ system \eqref{DNSE} 
has  a unique invariant measure $\mu^\nu$. Moreover
it is ergodic, i.e.
\begin{equation}\label{ergodic}
 \lim_{T\rightarrow\infty}\frac 1T\int_0^T \mathbb \varphi(\xi^\nu(t))dt=
 \int \varphi \ d\mu^{\nu} \qquad \text{ in } L^2(\Omega) 
\end{equation}
for any $\varphi \in \mathcal{C}_b(L^p)$ 
and initial vorticity in $ L^p$. 
Finally 
\begin{equation}\label{bal-mu} 
 \nu (p-1)\int \| \ |\xi|^{\frac p2 -1}\nabla \xi\|^2_{L^2} d\mu^\nu(\xi)
 +\gamma\int \|\xi\|^p_{L^p} d\mu^\nu(\xi) 
 = (p-1) Q \int \|\xi\|^{p-2}_{L^{p-2}} d\mu^\nu(\xi) .
\end{equation}
\end{theorem}
The latter equality comes from \eqref{Lp-formula}. Notice that this invariant measure $\mu^\nu$ is independent of $p$, since the assumption on the noise is independent of $p$.

\begin{remark} 
i) All the previous results hold true when 
$D$ is a smooth bounded domain in $ \mathbb R^2$, under  
the slip boundary condition coupled 
with a null vorticity on the boundary. In that case, the assumption 
on the noise has to be modified as
$\sum_{k\in \mathbb Z^2_0} |k|^2 q_k\|e_k\|^2_{L^p}<\infty$.
\\ii) For other conditions granting the uniqueness of the invariant measure see
e.g. \cite{BKL,DPZ,Deb,Ferrario99,Fla,FM94, HM11, KS02, M99}.
Anyway, our results hold when the noise is such that  the 
evolution of system \eqref{DNSE} is well defined for initial 
vorticity in $L^p$. In this case, we have that 
\eqref{bal-mu} is meaningful. 

In the following we shall fix $p= 4$; this allows to choose 
any kind of finite dimensional noise, whereas in the infinite dimensional case
($q_k \neq 0$ for all $k$) this is not a strong restriction. 
\end{remark}

Now, we fix the family of the unique invariant measures, as given in
Theorem \ref{teo-mu-unica}, and consider the limit of vanishing viscosity.

\begin{corollary} \label{ilcor}
Let $\gamma>0$. Then the family of  
invariant measures $\{\mu^{\nu}\}_{\nu>0}$ is tight in $H^{-s}$ for
any $s>0$;
 in particular there exists a measure $\mu^0$ in $H^{-s}$ such that
$$
 \mu^{\nu}\longrightarrow \mu^0 
 \quad {\rm weakly \ in }\quad H^{-s}
$$
as $\nu \longrightarrow 0$.
\end{corollary} 
\proof From \eqref{bal-mu} with $p=2$  we have
$$
 \int \|\xi\|^2_{H^0} d\mu^\nu(\xi) \le \frac {Q}{\gamma}
$$
uniformly in $\nu\in (0,\infty)$.
Then, using that $H^0$ is compactly
embedded in $H^{-s}$ we get tightness by means of  the Chebyshev 
inequality. \hfill $\Box$

\section{The  vanishing viscosity limit}
When $\nu=0$, we deal with the stochastic damped Euler equations
\begin{equation}\label{DEE}
\begin{cases}
d\xi^0+[u^0\cdot \nabla\xi^0 +\gamma \xi^0 ]\ dt=dw^{curl}&\\
\xi^0=\nabla^\perp \cdot u^0&
\end{cases}
\end{equation}
with periodic boundary conditions, as before. We always consider $\gamma>0$.

We are going to prove that this system has a stationary solution 
whose marginal at fixed time is the measure $\mu^0$ and
 that the following balance equation
holds:
$$
 \gamma \int \| \xi\|^2_{H^0}d\mu^0(\xi)=Q;
$$
moreover, considering the limit in the balance equation \eqref{bal-mu}
with $p=2$
we prove that
$$
 \lim_{\nu \to 0} \nu \int \|\nabla \xi\|^2_{H^0}d\mu^\nu(\xi)=0.
$$
This means that in the limit of vanishing viscosity, 
the damped stochastic equations \eqref{DNSE} have no dissipation of enstrophy.

However, instead of dealing with invariant measures, we  deal with
stationary processes (see next Remark \ref{oss-inv}).
Heuristically, we expect that there exists a stationary solution for the 
stochastic damped Euler system \eqref{DEE}, due to a balance between the
energy injected by the noise term and the dissipation of the
damping term.
More rigorously, in \cite{Bes2008}
it has  been  shown that the damped Euler equation  with 
a multiplicative noise has a stationary solution; there, the 
crucial estimate \eqref{stima22} was used that holds  for $\gamma>0$
(and $\nu \ge 0$).
The proof is even
easier with an additive noise; indeed, estimate \eqref{stima22} on the
finite dimensional approximating Galerkin system gives the existence
of an invariant measure by means of 
Krylov-Bogoliubov technique and 
we recover the existence of a
stationary solution for \eqref{DEE}.

Here, we want to investigate the properties for vanishing viscosity; in
particular the limit in the balance equation \eqref{bal-mu} with
$p=2$, that is
\begin{equation}\label{bal-mu-2}
 \nu\int \|\nabla \xi\|^2_{H^0} d\mu^\nu(\xi)
 +\gamma\int \|\xi\|^2_{H^0} d\mu^\nu(\xi) = Q
\end{equation}

Keeping in mind Corollary \ref{ilcor},
we consider the stationary stochastic process
$\xi^\nu$ whose law at any fixed 
time  is the measure $\mu^\nu$ of Theorem \ref{teo-mu-unica}, and 
take the limit of vanishing viscosity. We have

\begin{proposition}\label{pro-tight}
Let $s>0$.
The sequence $\{\xi^\nu\}_{\nu >0}$ of stationary processes solving 
\eqref{DNSE} has a subsequence
converging, as $\nu \to 0$,
 in $L^2_{loc}(0,\infty;H^{-s})\cap C([0,\infty);H^{-2-2s})$ 
(a.s.) to a process, which solves the damped
Euler system \eqref{DEE}. Moreover,  for any $p\ge 2$
the paths of the limit process
belong (a.s.) to $C([0,\infty); L^p_w)
\cap L^\infty_{loc}(0,\infty; L^p)$,
and the limit process is a stationary process in $L^p$. The marginal
at any fixed time of this limit process is the measure $\mu^0$.
\end{proposition}
\proof The proof is based on two steps: first we show that the 
sequence of the laws of $\xi^\nu$, $\nu >0$, is tight; 
then we pass to the limit in a suitable way and get that
the limit process is a weak solution of system \eqref{DEE}.
Notice that we find a weak solution to system \eqref{DEE} 
(in the probabilistic sense),
  whereas system \eqref{DNSE} has a unique strong solution.

Actually, the tightness and  the convergence of the stationary 
processes have already been done in \cite{Bes2008} 
for the damped Navier-Stokes equations with a multiplicative noise;
but there the analysis involved the velocity instead of the vorticity.
For the reader's convenience we recall the basic steps of the proof; 
the details can be found  in   \cite{Bes,Bes2008}.

Writing equation \eqref{DNSE} in the integral form 
$$
 \xi^\nu(t)=\xi^\nu(0)+\nu \int_0^t\Delta \xi^\nu(s)ds
 -\int_0^t u^\nu(s)\cdot \nabla\xi^\nu(s)ds
 -\gamma\int_0^t \xi^\nu(s)ds +w^{curl}(t),
$$
by usual estimations and bearing in mind estimate 
$$
 \sup_{0\le t<\infty}\mathbb E  \| \xi^\nu(t)\|_{L^4}^4\le C(4)
$$
from Theorem \ref{theorem1} (so we estimate 
$\sup_{0\le t<\infty}\mathbb E  \| \xi^\nu(t)\|_{H^0}^4$),
one gets that there exist constants $C$ and $C(p)$ such that
\[
\begin{split}
&\mathbb E\|\int_0^\cdot \Delta  \xi^\nu(s)ds\|^2_{W^{1,2}(0,T;H^{-2})}\le C\\
&\mathbb E\|\int_0^\cdot u^\nu(s)\cdot \nabla\xi^\nu(s)ds\|
   ^2_{W^{1,2}(0,T;H^{-2-s})}\le C \quad \text{ by
  }\eqref{cambiasegno} \text{ and }\eqref{stima-tril-xi0}\\
&\mathbb E\|\int_0^\cdot  \xi^\nu(s)ds\|^2_{W^{1,2}(0,T;H^0)}\le C\\
&\mathbb E\|w^{curl}\|^p_{W^{\alpha,p}(0,T;H^0)}\le C(p) 
     \quad \text{ by }\eqref{noise2}
\end{split}
\]
for some (and all) $s>0$,  $\alpha\in (0,\frac 12)$ and $p\ge 2$.
Therefore
$$
 \sup_{\nu\in (0,1)} \mathbb E \|\xi^\nu\|^2_{W^{\alpha,2}(0,T;H^{-2-s})}<\infty.
$$
On the other hand, we already know from  Theorem \ref{theorem1}  that
$$
 \sup_{\nu>0}  \mathbb E \|\xi^\nu\|^2_{L^2(0,T;H^0)}<\infty.
$$
Using that the space $L^2(0,T;H^0)\cap W^{\alpha,2}(0,T;H^{-2-s})$ 
is compactly embedded in $L^2(0,T;H^{-s})$ 
 (see, e.g., \cite{FlaGa}),
it follows that the sequence of laws of processes $\xi^\nu$ ($0<\nu<1$)
is tight in $L^2(0,T;H^{-s})$. On the other hand, using that both 
the spaces $ W^{1,2}(0,T;H^{-2-s})$ and $ W^{\alpha,p}(0,T;H^{-2-s})$
with $\alpha p>1$ 
are compactly embedded in   $C([0,T];H^{-2-2s})$,
we get tighness in $C([0,T];H^{-2-2s})$.

Let us endow $L^2_{loc}(0,\infty;H^{-s})$ by the distance
$$
 d_2(\xi,\zeta)=\sum_{n=1}^\infty 2^{-n} \min (\|\xi-\zeta\|_{L^2(0,n;H^{-s})},1)
$$
and $C([0,\infty);H^{-2-2s})$ by the distance
$$
 d_\infty(\xi,\zeta)=\sum_{n=1}^\infty 2^{-n} 
 \min (\|\xi-\zeta\|_{C([0,n];H^{-2-2s})},1).
$$
We have that the sequence $\{\xi^\nu\}$ is tight in
$L^2_{loc}(0,\infty;H^{-s})\cap C([0,\infty);H^{-2-2s})$.

From Prokhorov and Skorohod theorems  follows that 
there exists a basis $(\tilde\Omega, \tilde  F, \tilde P)$ 
and on this basis, $L^{2}_{loc}(0,\infty; H^{-s})
\cap C([0,\infty);H^{-2-2s})$-valued random variables 
$\tilde \xi^0$, $\tilde \xi^\nu$, such that
${\mathcal L}({\tilde \xi}^\nu)={\mathcal L}(\xi^\nu)$  on
$L^2_{loc}(0,\infty; H^{-s})\cap C([0,\infty);H^{-2-2s})$, and
\begin{equation}\label{conv-qo}
 \lim_{n \to \infty} {\tilde \xi}^{\nu_n} = {\tilde \xi}^0
 \qquad\text{ in } L^{2}_{loc}(0,\infty; H^{-s})
   \cap C([0,\infty);H^{-2-2s}),\; \tilde P-a.s. 
\end{equation}
for a subsequence with $\lim_{n \to \infty} \nu_n=0$.

The fact that the process $\tilde\xi^0$ solves system \eqref{DEE} is 
 classical. Indeed, considering $s=\frac 12$ we have that $\tilde \xi
 ^\nu \to \tilde \xi ^0$ in  
$L^{2}_{loc}(0,\infty; H^{-1/2})$; this means, according to \eqref{ell}, that 
$\tilde u^\nu \to \tilde u^0$ in $L^{2}_{loc}(0,\infty; H^{1/2})$. Since 
$H^{1/2}(D)\subset L^4(D)$, we get by estimates similar to \eqref{stima-tril-u}
that the quadratic term $[\tilde u^\nu \cdot \nabla]\tilde u^\nu$ 
converges weakly to 
$[\tilde u^0 \cdot \nabla]\tilde u^0$, i.e.
$$
 \int_D \int_0^t [\tilde u^\nu \cdot \nabla]\tilde u^\nu \cdot \psi \ ds\  dx
 \longrightarrow 
 \int_D \int_0^t [\tilde u^0 \cdot \nabla]\tilde u^0 \cdot \psi \ ds \ dx
 \qquad \tilde P-a.s.
$$
for all $t$ finite and $\psi \in [H^1]^2$.
For this it is enough to write
\begin{multline*}
\int_D \{[\tilde u^\nu \cdot \nabla]\tilde u^\nu \cdot \psi-
[\tilde u^0 \cdot \nabla]\tilde u^0 \cdot \psi\} \ dx
\\=
+ \int_D [(\tilde u^\nu -\tilde u^0)\cdot \nabla]\tilde u^\nu\cdot \psi \ dx
+ \int_D [\tilde u^0\cdot \nabla](\tilde u^\nu -\tilde u^0)  \cdot \psi\ dx.
\end{multline*}

In addition, $\tilde \xi^\nu$ and $\xi^\nu$ have the same law; then 
$\tilde \xi^\nu$ is a stationary process. By the convergence 
$\tilde P$-a.s. in $C([0,\infty);H^{-2-2s})$ 
we get that also $\tilde \xi^0$ is 
 a stationary process in $H^{-2-2s}$.

Finally,  from \eqref{stima-nu} we have that for $2\le p<\infty$
$$
 {\tilde \xi}^0 \in  L^\infty_{loc}(0,\infty;L^p) \qquad \tilde P-a.s. 
$$
Then, for $T<\infty$ almost each path $\tilde \xi^0 \in C([0,T];H^{-2-2s})\cap
L^{\infty}(0,T; L^p)$; thus it is weakly continuous in $L^p$,
i.e. we have  for any $\phi \in L^{p^\prime}$ ($\frac 1p +\frac 1{p^\prime}=1$) 
$$
 \lim_{t\to t_0} \int_D\tilde\xi^0(t)\phi\ dx = 
 \int_D \tilde\xi^0(t_0) \phi\ dx \qquad \tilde P-a.s.
$$ 
and for any $t \in [0,T]$
$$
 \|\tilde\xi^0(t)\|_{L^p} \le  \|\tilde\xi^0\|_{L^\infty(0,T;L^p)}
 \qquad \tilde P-a.s.
$$
(see \cite{Tem} p 263).\\
Hence, for every $t\ge 0$, the mapping 
$\tilde\omega\mapsto \tilde \xi^0(t,\tilde\omega)$ is well
defined from $\tilde\Omega$ to $L^p$ and it is weakly measurable. 
Since $L^p$
is a separable Banach space, it is strongly measurable (see
\cite{yosida} p 131). Therefore, it is meaningful to speak about the
law of $\tilde\xi^0(t)$ in $L^p$. 
The stationarity of $\tilde \xi^0$ in $L^p$ 
has to be understood in this sense.

By taking suitable subsequences we have
that $\mu^0$ is the law of $\tilde \xi^0(t)$ for any time $t$.
\hfill $\Box$

Let us denote by $\tilde \xi^0$ the stationary process 
solving \eqref{DEE}, as given in
Proposition \ref{pro-tight}. 
We have
\begin{proposition}\label{pro-bil}
For  any time $t$
\begin{equation}\label{bilancioE}
 \gamma \tilde{\mathbb  E}\| \tilde\xi^0(t)\|^2_{H^0}=Q.
\end{equation}
\end{proposition}
\proof
Choosing $p=4$ in \eqref{stima22} of Theorem \ref{theorem1} we have
$$
  \tilde {\mathbb E} \|\tilde \xi^\nu(t)\|^4_{L^4}\le C(4) .
$$
This bound implies  
$$
 {\tilde \xi}^\nu(t) \longrightarrow \tilde \xi^0(t) \qquad 
 \text{weakly in  }L^4(\tilde \Omega \times D);
$$
for the limit we have
\begin{equation}\label{stima2}
 \tilde {\mathbb E} \|\tilde \xi^0(t)\|^4_{L^4}\le 
  \liminf_{\nu \to 0} \tilde {\mathbb E} \|\tilde \xi^\nu(t)\|^4_{L^4}\le C(4).
\end{equation} 
By working on the first equation of \eqref{DEE}, It\^o formula 
for $d\|\tilde \xi^0(t)\|^2_{H^0}$ provides
\begin{equation}
 \|\tilde\xi^0(t)\|_{H^0}^2 +\gamma \int_0^t \|\tilde\xi^0(s)\|_{H^0}^2 ds
 =\|\tilde\xi^0(0)\|^2_{H^0}
 +tQ +\int_0^t \langle\tilde\xi^0(s),d\tilde w^{curl}(s)\rangle,
\end{equation}
$\tilde P$-a.s.
For this we have used \eqref{cons}, having that, for any $s$,
$\tilde\xi^0(s) \in L^4(D)$ a.s. from \eqref{stima2}.
\\
Taking expectation and using stationarity we get \eqref{bilancioE}.
\hfill$\Box$

Equation \eqref{bilancioE} can be rewritten as
$$
 \gamma \int \|\xi\|_{H^0}^2 d\mu^0(\xi)=Q.
$$
\begin{remark}\label{oss-inv}
At this point, we are not able to prove that $\mu^{0}$ is an 
invariant measure for the system \eqref{DEE}. In fact, the 
transition semigroup associated to \eqref{DEE} can not 
be defined in $H^{0}$:   
existence of a solution holds for initial vorticity in
$H^0$ but uniqueness requires stronger assumptions (see
\cite{Bes2000} and \cite{BesFla2000}). But to get the Feller and Markov properties
in a space smaller than $H^0$ is not trivial. Some work in progress 
in that direction is being made by the current authors.  
\end{remark}

Now we have our main result
\begin{theorem}
For any $\gamma>0$, we have
\begin{equation}\label{tesi}
\lim_{\nu \to 0}\; \nu \int \|\nabla \xi\|^2_{H^0}d\mu^\nu(\xi)=0.
\end{equation}
\end{theorem}
\proof
Let us write the balance equation \eqref{bal-mu-2} 
in terms of the stationary process $\xi^\nu$, at any fixed time $t$:
\begin{equation}\label{bal}
 \nu\mathbb E \| \nabla \xi^\nu(t)\|^2_{H^0} 
 +\gamma\mathbb E \|\xi^\nu(t)\|^2_{H^0}  =
 Q.
\end{equation}
Considering the weak limit as in Proposition \ref{pro-tight} 
and \ref{pro-bil} we have
\begin{equation}
\begin{split}
 \limsup_{\nu \to 0}\;\nu \tilde{\mathbb E}
    \|\nabla \tilde\xi^\nu(t)\|^2_{H^0}
&=
 Q
 -\gamma\liminf_{\nu \to 0} \tilde{\mathbb E} \|\tilde\xi^\nu(t)\|^2_{H^0}
\\&
\le Q-\gamma \tilde{\mathbb E} \|\tilde\xi^0(t)\|^2_{H^0}
 \text{ by }\eqref{stima2}
\\&
 =0 \text{ by } \eqref{bilancioE}.
\end{split}
\end{equation}
This gives \eqref{tesi}. \hfill $\Box$

From this result we obtain the convergence of the mean  enstrophy.
\begin{corollary}
For any $\gamma>0$, we have
\begin{equation}
 \lim_{\nu \to 0} \int \|\xi\|^2_{H^0}d\mu^\nu(\xi)
 = 
 \int \|\xi\|^2_{H^0}d\mu^0(\xi).
\end{equation}
\end{corollary}
\proof
We consider the limit as $\nu \to 0$ in \eqref{bal-mu-2}; 
then use \eqref{tesi} and \eqref{bilancioE}. \hfill $\Box$

\begin{remark}
All the results proved for the enstrophy $\xi$ can be repeated and hence
hold for the velocity $u$; norms of one order less of regularity are 
involved and therefore the proofs are even easier.  
This means in particular that for the stochastic  damped 2D Navier-Stokes  
equations,  there is no anomalous dissipation of energy as $\nu \to 0$
and energy balance equation holds for $\nu>0$ and also $\nu=0$.
\end{remark}
{\bf Acknowledgment:} The work of H. Bessaih was supported in part by the
GNAMPA-INDAM project "Professori Visitatori". We would like to thank
the hospitality of the Department of Mathematics of the University of Pavia
where part of this research started and the IMA in Minneapolis where 
the paper has been finalized.

\end{document}